\documentclass{amsart}
\newtheorem{theorem}{Theorem}[section]

\newtheorem{lemma}[theorem]{Lemma}

\theoremstyle{definition}
\newtheorem{definition}[theorem]{Definition}
\theoremstyle{remark}
\newtheorem{remark}[theorem]{Remark}
\theoremstyle{example}

\begin{document}

\parskip0.2cm

\title{Nambu structures and integrable 1-forms}

\author{Jean-Paul Dufour and Mikhail Zhitomirskii \\
}

\date{october 1999}

\newpage
\maketitle

\begin{abstract} {\normalsize
Some years ago Mosh\'e Flato pointed out that it could be interesting to develop
the Nambu's idea (\cite{Nambu}) to generalize Hamiltonian mechanic.
An interesting new formalism in that direction
was proposed by L. Takhtajan (\cite{T}).
 His  theory gave new perspectives concerning deformation quantization,
 and many
authors have developed its mathematical features (\cite{D-F}, \cite{D-F-S-T}).

The purpose of this paper is to show that this theory,
at first dedicated to physic,
gives a new point of view for the
study of singularities of
 integrable 1-forms.

Namely, we will prove that any
integrable 1-form which vanishes at a point
and has a non-zero linear part at this point is,
up to multiplication by a non-vanishing function, the
formal pull-back of a two dimensional 1-form.
We also obtain a classification of quadratic
integrable 1-forms.}

\end{abstract}




\vskip1cm
{\sl To the memory of Mosh\'e Flato}

\vskip1cm

\section{Generalities}

In the reference \cite{T}, L. Takhtajan,  in 1994,  proposed
a formalism which
generalizes the Poisson bracket. Let
$M$ be a manifold and $A$ the algebra of smooth functions on
$M$. A {\bf Nambu structure of order} $r$ on  $M$ is an $r$-linear
skew-symmetric
map
$$ A\times\cdots \times A\rightarrow A :$$
$$(f_1,\dots ,f_r)\mapsto \{f_1,\dots ,f_r\} $$ which satisfies
the following properties:
$$\{ f_1,\dots ,f_{r-1},gh\} = \{ f_1,\dots ,f_{r-1},g\} h+
g\{ f_1,\dots ,f_{r-1},h\}\eqno{(L)}$$
$$\{ f_1,\dots ,f_{r-1},\{ g_1,\dots ,g_{r}\}\}=
$$$$\sum_{i=1}^r \{ g_1,\dots ,g_{i-1},\{ f_1,\dots ,f_{r-1},g_i\} ,g_{i+1},\dots
,g_{r}\}\eqno{(FI)}$$
for any $f_1,$...,$f_{r-1},$ $g,$ $h,$ $g_1,$...,$g_r$ in $A.$

In this definition $(L)$ stands for \it Leibniz property\rm , $(FI)$
for \it fundamental identity \rm  or for \it Filippov's identity \rm
(see \cite{M-V-V}).
For $r=2,$ $(FI)$ is just  Jacobi's identity, so a Nambu
structure of order 2 is a Poisson structure.

The identity $(L)$ implies that $X_{f_1\cdots f_{r-1}}:g
\mapsto \{ f_1,\dots ,f_{r-1},g\}$ is a derivation
of $A,$ hence a vector field on $M$: It is,
by definition, the {\it Hamiltonian vector field}
associated to $f_1,\dots ,f_{r-1}.$

The identity $(L)$ also implies
that there is an $r$-vector field $\Lambda$  such that
$$\{ f_1,\dots ,f_{r}\} =\Lambda (df_1,\dots ,df_r).$$
This $\Lambda$ is called a {\bf Nambu tensor.}
We can also consider the usual vector fields
as Nambu structures
of order 1.

The identity $(FI)$ implies that  Hamiltonian vector
fields define an integrable distribution,
like in Poisson's case. So, we have on $M$ a singular foliation which
generalizes symplectic foliations of Poisson manifolds.

Since 1996 appeared three  proofs of the following
surprising result (\cite{G}, \cite{A-G},
\cite{N}).

\begin{theorem} [Local Triviality Theorem]
Let $\Lambda$ be a Nambu tensor of order $r>2.$ Near any point
at which $\Lambda$ does not vanish there  are local coordinates
$x_1,\dots ,x_n$ such that
$$\Lambda ={\partial \over\partial x_1}\wedge \dots
\wedge {\partial \over\partial x_r}.$$
\end{theorem}

In particular this theorem
shows that there are only two types of leaf for the foliation associated
to $\Lambda$:
Either it reduces to a point (zero of $\Lambda$) or it is $r$-dimensional.

This theorem leads to a ``covariant'' presentation of Nambu tensors.
Suppose that we have a volume form $\Omega$ on
our manifold $M.$  Set $\omega:=i_\Lambda \Omega .$ Then we have the
following result (\cite{D-N}).

\begin{theorem} Suppose $\Lambda$ is a $r$-vector on $M$ such
that either $r>2$ or $r=2$ but, in this case, maximal rank of $\Lambda$ is 2.
If $r$ is equal to the dimension $n$ of $M$,
then $\Lambda$ is always a Nambu tensor.
When $r<n,$ $\Lambda$ is a Nambu tensor if and only  if we have
$$i_A\omega\wedge\omega =0$$
$$i_A\omega\wedge d\omega =0$$
for every $(n-r-1)$-vector $A.$
\end{theorem}

The first relation in this theorem says that $\omega$ is decomposable
at each point, the second is an ``integrability'' condition. In the case
$r=n-1,$ $\omega$ is just an integrable 1-form, i.e., a 1-form
such that $\omega\wedge d\omega=0.$ In the case
$r<n-1,$ $\omega$ can be called an {\sl integrable $(n-r)$-form}, see
\cite{Me}.
Roughly speaking, this theorem  says that a
Nambu structure (or a Poisson structure of maximal rank 2) is
exactly the ``dual'' of  an integrable $p-$form.

For Nambu structures there is an analogous of the so called
modular vector field
(\cite{W}, \cite{D-H}) which can be defined as follows.

\begin{definition}
Let $\Lambda$ be a Nambu tensor of order $r$ and $\Omega$ be
a volume form on the manifold $M.$
The {\bf modular tensor} of $\Lambda$ with respect to $\Omega$
is the tensor field $D_\Omega\Lambda$
defined by the formula
$$i_{D_\Omega\Lambda}\Omega =d(i_\Lambda\Omega ).$$
\end{definition}

Using the local triviality theorem we can prove the following results.

\begin{theorem} \label{proprietesducurl} The modular tensors
are also Nambu tensors. If $\Lambda$ is a Nambu
tensor of order $r$ with $r>2$ or with $r=2$, but with maximal rank 2,
then we have, for any
volume form $\Omega $, for every $s,$ $s=0,$ 1,..., $r-2,$
and for any smooth
functions  $g_1, \dots , g_s$, the following properties

\noindent $1)\ i_{(dg_1 \wedge \dots \wedge dg_s)}
D_{\Omega}\Lambda\wedge\Lambda =0,$

\noindent $ 2)\ [i_{(dg_1 \wedge \dots \wedge dg_s)}
D_{\Omega}\Lambda , \Lambda ]=0,$
where the bracket $[\ ,\ ]$, is the Schouten bracket.
\end{theorem}

Note that the property 2)(with $s=0$) remains valid for
any Poisson tensor, even if its maximal rank
is more than 2.

\section{The Kupka phenomenon}

The Kupka phenomenon (\cite{Kupka})
is the following: If $\omega$ is an integrable 1-form
such that $d\omega$ is non zero at a point, then near this point
there are local coordinates
$x_1,\dots ,x_n$ such that $\omega$ depends
 only on two variables, i.e., we have
$$\omega =a(x_1,x_2)dx_1+b(x_1,x_2)dx_2.$$

 Using the fact that integrable 1-forms are the ``duals'' of
 Nambu tensors of
order $n-1$ ($n$ is the dimension of the ambiant manifold), we
could rewrite this result in terms of Nambu tensors,
but, hereafter, we will give a generalization of this result. For this we will
use the following vocabulary.

\begin{definition} \label{2.r} Let $A$ be a Nambu tensor.
We will say that $A$ is of
type $2.r$ if there are $r$ commuting and everywhere linearly independent
vector fields $X_1,\dots ,X_r$ such that we have
$$X_i\wedge A=0$$$$[X_i,A]=0$$
for every $i=1,\dots ,r ($[\ ,\ ]$ is the Schouten bracket).$
\end{definition}

\begin{remark} Locally this means that there are local
coordinates $x_1,\dots ,x_n$
such that
$$A={\partial /\partial x_1}\wedge \dots\wedge{\partial /\partial x_r}\wedge B$$
where $B$ is a Nambu tensor independent of the coordinates $x_1,\dots ,x_r.$
\end{remark}

\begin{theorem} [generalized Kupka phenomenon] \label{KupkaNambu}
Let $\Lambda$ be a
Nambu tensor and $\Omega$ a volume form. If  $D_\Omega\Lambda$ is a.e.
non zero
and is of type $2.r$ in a neighborhood of a point $m$
then $\Lambda$ is also of type
$2.r$ in a (possibly different) neighborhood of $m.$
\end{theorem}

{\sl Proof.} We can choose local coordinates $(x_1,\dots ,x_n)$ such that
$X_i=\partial /\partial x_i$ for $i=1,\dots ,r$ and $\Omega =dx_1\wedge \dots\wedge dx_n.$
Then we have
$$D_\Omega\Lambda =
{\partial
/\partial x_1}\wedge \dots\wedge{\partial /\partial x_r}\wedge Y$$
where
$$Y=\sum Y_{i_1\dots i_{q-1-r}}{\partial /\partial x_{i_1}}\wedge \dots\wedge
{\partial /\partial x_{i_{q-1-r}}}$$
is a $(q-1-r)$-tensor field independent of $x_1,$...,$x_r.$
Since $D_\Omega\Lambda$
is a.e. non zero we can suppose that $Z:=Y_{(r+1)\dots (q-1)}$ is a.e. non zero.

Set $\nu =dx_1\wedge \dots\wedge dx_{i-1}\wedge dx_{i+1}
\wedge\dots\wedge dx_{q-1}.$ We have
$i_\nu (D_\Omega\Lambda)=Z\partial /\partial x_i.$ The relation 1) of theorem
 \ref{proprietesducurl}, implies $\partial /\partial x_i\wedge \Lambda=0 .$
The latter relation  holds for $i=1,\dots ,r,$ so we obtain
$$\Lambda =
{\partial /\partial x_1}\wedge \dots\wedge{\partial /\partial x_r}\wedge P,$$
where $P$ is a $(q-r)$-tensor field.

Since $Z$ is independent of $x_1,\dots ,x_r,$ the relation 2) of
theorem \ref{proprietesducurl}
implies that  $[(\partial /\partial x_i, P ]=0 .$ It follows that
$P$ is independent of $x_1,\dots ,x_r.$ This ends the proof of our theorem.
\hfill $\triangle$

Let $\Lambda$ be the Nambu tensor of order $n-1$ associated with
an integrable 1-form $\omega ,$ such that
$d\omega\neq 0$ at a point $m.$
Then the modular tensor of $\Lambda$ is non-zero at $m$
and the local triviality theorem for
regular Nambu structures says that it is locally of the form
${\partial /\partial x_1}\wedge \dots\wedge{\partial /\partial x_{n-2}},$
so it is of type $2.(n-2).$
The theorem above says that $\Lambda$ is also of type $2.(n-2).$
According to the preceding remark we have, locally,
$$\Lambda={\partial /\partial x_1}\wedge \dots\wedge{\partial /\partial x_{n-2}}\wedge B$$
where $B$ is  independent of the coordinates $x_1,\dots ,x_{n-2}.$
Therefore,
up to multiplication
by a non-vanishing function, $\omega$ depends only on $2$ coordinates.
It is easy to see that the latter remains true without
multiplication by a nonvanishing function under
a suitable choice of the involved
volume form.  Therefore
our theorem
can be thought as a generalization of the Kupka phenomenon.

For example, the formulated theorem has the following
corollary (which can be proved directly).

\begin{theorem} Let $\omega$ be an integrable 1-form on ${\mathbb R}^n$ or
${\mathbb C}^n$. If $d\omega$
is a.e. non zero and depends on less than $s$ coordinates
in a neighborhood of 0 then
 we have the same for $\omega .$
\end{theorem}

\section{Nambu tensors of order $n-1$ with a non-zero linear part}

In this section we give a formal normal form for Nambu tensors of
order $n-1$, vanishing
at a point $m$, but with a non-zero linear part at that point;
it generalizes the one we gave
in \cite{D-Z} for the 3 dimensional case.

We will distinguish the following two cases.

The simple case is the one where the modular tensor doesn't vanish:
In that case by  ``Kupka phenomenon'' our Nambu tensor has
the local form
$${\partial /\partial x_1}\wedge \dots\wedge{\partial /\partial x_{n-2}}\wedge X$$
where $X$ is a vector field independent of the coordinates
$x_1,\dots ,x_{n-1}.$ Thus the
local classification of Nambu tensors reduces to that of
2-dimensional vector fields, up to orbital equivalence.

The difficult case is the one where the modular tensor
vanishes at $m.$ In this case we have the following
theorem.

\begin{theorem}\label{NambuNF} Let $\Lambda$ be
a Nambu tensor of order $n-1$ on an $n$-dimensional manifold with
$n\geq 3.$ Suppose that $\Lambda$ vanishes at a point $m$,  but
has a non-zero linear part at this point. Suppose also that the
modular tensors of $\Lambda $ vanish at $m$. Then there are local
coordinates $x_1,\dots ,x_n,$ in a neighborhood of $m$ such that
$$\{x_1,\dots ,x_{n-1}\}=x_n$$ $$\{x_1,\dots
,x_{i-1},x_{i+1},\dots ,x_n\}=(-1)^{n-i}({\partial f /\partial
x_i}+ x_n{\partial g /\partial x_i})+\epsilon_i,$$ for $i=1,\dots
,n-1,$ where $f$ and $g$ are smooth functions, independent of
$x_n,$
 such that $df\wedge dg=0$,
and $\epsilon_i$ is a smooth flat function at the origin
(i.e., his Taylor expansion
vanishes at $m.$)
\end{theorem}

The sequel of this section is dedicated to the proof of this theorem.

\noindent {\sl Study of the linear part of $\Lambda .$}
According to \cite{D-N}
the linear part $\Lambda^{(1)}$ has, in a suitable coordinates
system, one of the
following normal forms.

\noindent{\it Type 1}:
$$ \Lambda^{(1)}=\sum_{i=1}^{r}\pm x_i{\partial  /\partial x_1}\wedge\cdots\wedge
 {\partial /\partial x_{i-1}}\wedge {\partial /\partial x_{i+1}}\wedge\cdots
\wedge {\partial /\partial x_n},$$
which corresponds to a linear integrable 1-form of type $d( \sum_{i=1}^{r}\pm x_i^2/2).$

\noindent{\it Type 2}: $$ \Lambda^{(1)}=\partial /\partial
x_1\wedge\cdots\wedge
  {\partial /\partial x_{n-2}}\wedge X^{(1)},$$
where  $X^{(1)}$ is a zero-trace linear vector field depending only on $x_{n-1}$ and $x_n.$
This normal form  corresponds to a linear integrable 1-form
depending only on $x_{n-1}$ and $x_n.$

An elementary calculation shows that, in each of the cases, there are
(possibly) new linear
coordinates with
\begin{equation}\label{NF0}
\{x_1,\dots ,x_{n-1}\}^{(1)}=x_n
\end{equation}
for the linear Nambu structure determined by $\Lambda^{(1)}.$
This means that the
associated 1-form is of type $x_ndx_n+\sum_{i=1}^{n-1}l_idx_i.$

\begin{remark} In fact the preceding theorem is true for every case
 where one can find coordinates
satisfying (\ref{NF0}). The only case where this is not so
is the type 2 case with $X^{(1)}$
equivalent to $x_{n-1}{\partial /\partial x_{n-1}}+x_n{\partial /\partial x_{n}}.$
\end{remark}

In the sequel of the proof of theorem \ref{NambuNF}
we will use following notations:
$$x:=(x_1,\dots ,x_{n-1}),\ \ y:=x_n.$$
We also develop the function $h(x_1,\dots ,x_n)=:h(x,y)$ in the form
$$h^{(0)}+h^{(1)}+\cdots +h^{(p)}+\cdots ,$$
where $h^{(p)}$  is a $p$-homogeneous
polynomial in $x_1,\dots ,x_{n-1}$ with coefficients
depending smoothly on $y$ ($y$ is considered as a parameter).

\begin{lemma} \label{lemmaNF} Let $r\geq 0.$
Suppose that there are coordinates
$x =(x_1,\dots ,x_{n-1})$ and $y$ such that
$$\{x_1,\dots ,x_{n-1}\}=y+c^{(r+2)}(x,y)+c^{(r+3)}(x,y)+\cdots$$
$$\{x_1,\dots ,x_{i-1},x_{i+1},\dots ,x_{n-1},y\}=(-1)^{n-i}(a_i^{(0)}(x,y)
+a_i^{(1)}(x,y)+\cdots )$$
where $a_i^{(0)},\dots ,a_i^{(r-1)}$ are affine
with respect to $y$ (vacuous hypothesis for
$r=0$). There is a coordinates transformation of the form
$$x'_1=x_1+\mu^{(r+2)}(x,y)$$
$$x'_2=x_2,\dots ,x'_{n-1}=x_{n-1}$$
$$y'=y+\gamma^{(r+1)}(x,y)+\gamma^{(r+2)}(x,y)$$
which gives
$$\{x'_1,\dots ,x'_{n-1}\}=y'+C^{(r+3)}(x',y')+C^{(r+4)}(x',y')+\cdots$$
$$\{x'_1,\dots ,x'_{i-1},x'_{i+1},\dots ,x'_{n-1},y'\}=
(-1)^{n-i}(a_i^{(0)}(x',y')+\cdots +a_i^{(r-1)}(x',y')+$$$$
A_i^{(r)}(x',y')+A_i^{(r+1)}(x',y')\cdots )$$
where $A_i^{(r)}$ is affine in $y'.$
\end{lemma}

\noindent {\sl Proof of the lemma.} Make a coordinates transformation of the form
$\tilde x=x$, $\tilde y=y(1+e^{(r+1)}(x,y)).$ We obtain
$$\{\tilde x_1,\dots ,\tilde x_{n-1}\}=\tilde y+\tilde y(\tilde c^{(r+1)}(\tilde x,\tilde y))
+\tilde c^{(r+2)}(\tilde x,\tilde y)+\cdots$$
$$\{\tilde x_1,\dots ,\tilde x_{i-1},\tilde x_{i+1},\dots ,\tilde x_{n-1},\tilde y\}=
(-1)^{n-i}(a_i^{(0)}(\tilde x,\tilde y)
+\cdots +a_i^{(r-1)}(\tilde x,\tilde y)+$$$$
A_i^{(r)}(\tilde x,\tilde y)+\cdots )$$
with
\begin{equation}\label{Ar}
A_i^{(r)}=a_i^{(r)} -y^2{\partial e^{r+1}/\partial x_i}.
\end{equation}

Now denoting $\Omega =dx_1\wedge \dots\wedge dx_{n-1}\wedge dy,$ we have
$\omega :=i_\Lambda\Omega=\Gamma dy+\sum_i\Delta_idx_i$ with
$$\Gamma =\{x_1,\dots ,x_{n-1}\},\ \Delta_i=(-1)^{n-i}
\{x_1,\dots ,x_{i-1},x_{i+1},\dots ,x_n\}.$$
Recall that we have $\omega\wedge d\omega =0.$
The terms with $dx_i\wedge dx_j\wedge dy$
in this last equation give
\begin{equation}\label{Ar+}
\Gamma ({\partial \Delta_i /\partial x_j}-{\partial \Delta_j /\partial x_i})+
\Delta_i({\partial \Delta_j /\partial y}-{\partial \Gamma /\partial x_j})-
\Delta_j({\partial \Delta_i /\partial y}-{\partial \Gamma /\partial x_i})=0.
\end{equation}

Express $\Delta _k$ in the form
 $\Delta_k=\alpha_k(x)+y\beta_k(x)+y^2\delta_k(x,y)$. Our hypothesis
says that $\delta_k$ have developments $\delta_k^{(r)}+\delta_k^{(r+1)}+\cdots . $
Now, if we compare terms with $y^3$ and of order $r-1$ in the preceding
equation, we get
\begin{equation}\label{Ar++}
{\partial \delta_i^{(r)} /\partial x_j}-{\partial \delta_j^{(r)} /\partial x_i}
=0.
  \end{equation}
Equation (\ref{Ar}) can be rewritten in the form
$$A_i^{(r)}=\alpha_i^{(r)}+\beta_i^{(r)} +y^2(\delta_i^{(r)}-
{\partial e^{r+1} /\partial x_i}). $$ By the Poincar\'e lemma we
can choose $e^{(r+1)}$ such that $A_i^{(r)}$ are affine in $y$ (we
erase  $\delta_i^{(r)}$).

Now (after this  coordinates transformation) we can suppose
$$\Gamma =\{x_1,\dots ,x_{n-1}\}=y+yc^{(r+1)}(x,y)+c^{(r+2)}(x,y)+\cdots$$
$$\Delta_i=(-1)^{n-i}\{x_1,\dots ,x_{i-1},x_{i+1},\dots ,x_{n-1},y\}=a_i^{(0)}(x,y)
+a_i^{(1)}(x,y)+\cdots ,$$
where $a_i^{(s)}$ are affine in $y$ for $s=0,\dots ,r.$

In a second step we use a  coordinates transformation of the form $\tilde x_1 =
x_1+\theta^{(r+2)}(x,y),$ $\tilde x_2=x_2,\dots ,\tilde x_{n-1}=x_{n-1},\tilde y=y$
with $\partial\theta^{(r+2)}/\partial x_1=-c^{(r+1)}.$ Then we obtain
$$\{\tilde x_1,\dots ,\tilde x_{n-1}\}=\tilde y+
+\tilde c^{(r+2)}(\tilde x,\tilde y)+\cdots$$
$$\{\tilde x_1,\dots ,\tilde x_{i-1},\tilde x_{i+1},\dots ,\tilde x_{n-1},\tilde y\}=
(-1)^{n-i}(a_i^{(0)}(\tilde x,\tilde y)
+\cdots +a_i^{(r)}(\tilde x,\tilde y)+$$$$
A_i^{(r+1)}(\tilde x,\tilde y)+\cdots ).$$
Now we can suppose
$$\Gamma =y+c^{(r+2)}+\cdots$$
$$\Delta_i=a_i^{(0)}+\cdots ,$$
where the $a_i^{(s)}$ are affine in $y$ for $s=0,\dots ,r.$

Finally, to achieve the proof of the lemma,
it suffices to perform a  coordinates transformation
$\tilde x=x,$ $\tilde y =y+c^{(r+2)}. $
\hfill $\triangle$

We continue  the proof of theorem \ref{NambuNF}.

Since we have (\ref{NF0}), we can take $\{x_1,\dots ,x_{n-1}\}$ as a
new variable $y$ to get $$\{x_1,\dots ,x_{n-1}\}=y.$$ Then the
hypothesis of lemma \ref{lemmaNF} holds for $r=0.$ We can apply
inductively this lemma  to show
that, after a formal  coordinates transformation(the formal composition
of the coordinates transformations given by the lemma), we obtain
$$\{x_1,\dots ,x_{n-1}\}=y$$ $$\{x_1,\dots ,x_{i-1},x_{i+1},\dots
,x_{n-1},y\}=(-1)^{n-i}A_i, $$ where the functions $A_i$ have
formal developments $A_i^{(0)}+A_i^{(1)}+\cdots$ with  all terms
here being affine in $y.$ Therefore  we can suppose that we have
formally $$\{x_1,\dots ,x_{i-1},x_{i+1},\dots
,x_{n-1},y\}=(-1)^{n-i}(\alpha_i(x)+y\beta_i(x))$$ for
$i=1,\dots ,n-1.$

Set
$\Omega =dx_1\wedge \dots\wedge dx_{n-1}\wedge dy$. Then
the associated integrable 1-form $\omega $ has the form
$$\omega =\sum_{i=1}^{n-1}(\alpha_i(x)+y\beta_i(x))dx_i+ydy.$$
The equation $\omega\wedge d\omega =0$ implies
$$\alpha_i\beta_j-\alpha_j\beta_i\pm ({\partial \alpha_i /\partial x_j}-
{\partial \alpha_j /\partial x_i})\pm ({\partial \beta_i /\partial x_j}-
{\partial \beta_j /\partial x_i})=0.$$
So we obtain, for every $i$ and $j,$
$${\partial \alpha_i /\partial x_j}={\partial \alpha_j /\partial x_i},\
{\partial \beta_i /\partial x_j}={\partial \beta_j /\partial x_i},\
\alpha_i\beta_j=\alpha_j\beta_i.$$
The Poincar\'e lemma gives $\alpha_i=\partial f/\partial x_i,$
$\alpha_i=\partial g/\partial x_i,$ for every $i.$

Therefore the latter equations leads to
$df\wedge dg=0$. This ends the proof of theorem \ref{NambuNF}
\hfill $\triangle$

Theorem \ref{NambuNF} has the following
consequence concerning integrable 1-forms.

\begin{theorem} \label{1-formNF} Let $\omega$ be an integrable 1-form
which vanishes
at a point $m$ and has a non-zero linear part at this point.
Then, up to multiplication by a
non-vanishing function, $\omega$ is, formally, the pullback of an integrable
1-form depending only on 2 variables.
\end{theorem}

\noindent {\sl Proof.} If $d\omega$ is non-zero, we can
apply the Kupka phenomenon.
If $d\omega$ vanishes at $m$ then the Nambu vector associated to $\omega$ has
the formal form of theorem \ref{NambuNF}. So we can suppose that
$$\Lambda =y{\partial /\partial x_1}\wedge\cdots\wedge
  {\partial /\partial x_{n-1}}+(-1)^{n-i}({\partial f /\partial x_i}
+\sum y{\partial g /\partial x_i}){\partial /\partial x_1}\wedge\cdots
$$$$\cdots\wedge
  {\partial /\partial x_{i-1}}\wedge{\partial /\partial x_{i+1}}
\wedge\cdots\wedge
  {\partial /\partial x_{n-1}}\wedge{\partial /\partial y}.$$
Therefore
$$\omega =df +ydg+ydy$$
up to multiplication by a non-vanishing function
(the Jacobian of the change
 of coordinates).

Since we also have $df\wedge dg=0$ we can apply the result of
\cite{Moussu} to exhibit a function $h(x)$ such that
$$f=a\circ h,\ g=b\circ h$$ (at least at the level of
formal series; here $a$ and $b$ are
functions in one variable). Then we have
$$\omega =(a'(h)+yb'(h))dh +ydy =\phi^*\omega_2$$
with $\omega_2 =(a'(u)+vb'(u))du+vdv$ and $\phi :(x,y)\mapsto (h(x),y).$
This ends the proof of the theorem.
\hfill $\triangle$

\begin{remark} Theorems \ref{NambuNF} and \ref{1-formNF} give only {\sl formal}
normal forms for $(n-1)$ order Nambu structures or integrable 1-forms.
We do not know
if there are smooth or analytic versions.
\end{remark}

\begin{remark} In fact theorem \ref{1-formNF} can be proven directly (without
using Nambu formalism). The crucial point of the proof is that, up to
multiplication by a non-vanishing function, an
integrable 1-form $ydy +\sum_{i}A_idx_i$
is formally equivalent to a form $ydy +\alpha_0+y\alpha_1,$
where $\alpha_0$ and $\alpha_1$
are 1-forms depending on $x_1,$...,$x_{n-1}$ only. This
result has the following generalization.
\end{remark}

\begin{theorem} Let
$\omega =y^pdy +\sum_{i=1}^{n-1}A_idx_i$ be
an integrable 1-form on ${\mathbb R}^n$ (or
${\mathbb C}^n$). Then, up to multiplication by a non-vanishing function,
$\omega$
is formally equivalent to an integrable 1-form
$$\omega_0=y^pdy +\sum_{i=0}^{p}y^i\alpha_i$$
where $\alpha_i$ are 1-forms depending only on $x_1,$..., $x_{n-1}.$
\end{theorem}

{\sl Proof.} We consider the associated Nambu tensor
$$\Lambda =y^{p}{\partial /\partial x_1}\wedge\cdots\wedge
  {\partial /\partial x_{n-1}}+\sum (-1)^{n-i}A_i)
{\partial /\partial x_1}\wedge\cdots
\wedge
  {\partial /\partial x_{i-1}}\wedge{\partial /\partial x_{i+1}}
\wedge\cdots$$$$\cdots\wedge
  {\partial /\partial x_{n-1}}\wedge{\partial /\partial y}.$$
If $A_i^{(0)},\dots ,A_i^{(r-1)}$ are all polynomials
of degree $p$ in $y$, with coefficients
depending on $x$, then
we can apply exactly the same method as in the first step of
the proof of lemma \ref{lemmaNF} to bring  $A_i^{(r)}$
to a polynomial in $y$ of
degree $p$. In order to get this, we make
a coordinates transformation $\tilde x=x,$
$\tilde y=y(1+c^{r+1})$ with the notations of the proof of this lemma.
Then
$$\{\tilde x_1,\dots ,\tilde x_{n-1}\}=\tilde y^p
(1+\tilde c^{(r+1)}(\tilde x,\tilde y)+\cdots )^p$$
$$\{\tilde x_1,\dots ,\tilde x_{i-1},\tilde x_{i+1},\dots ,\tilde x_{n-1},\tilde y\}=
(-1)^{n-i}A_i(1+c^{(r+1)}-y{\partial c^{(r+1)} / \partial y})
-$$$$(-1)^{n-i}y^{p+1}{\partial c^{(r+1)}/ \partial x_i}
=(-1)^{n-i}(A_i^{(0)}(\tilde x,\tilde y)+\cdots +$$$$A_i^{(r-1)}(\tilde x,\tilde y)+
A_i^{(r-1)}(\tilde x,\tilde y)
-\tilde y^{p+1}{\partial c^{(r+1)}/ \partial x_i}+\tilde A_i^{(r+1)}+\cdots$$
Now we develop $A_i^{(r)}$ in the form
$$\alpha_{i,0}^{(r)}+y\alpha_{i,1}^{(r)}+\cdots +y^p\alpha_{i,p}^{(r)}+y^{p+1}\delta_i^{(r)}$$
where $\alpha_{i,j}^{(r)}$ depends only on $x$ for $j=0,\dots ,p.$

The identity $\omega\wedge d\omega =0$ implies that
$$y^p( {\partial A_i/ \partial x_j}-{\partial A_j/ \partial x_i})+
A_i{\partial A_j/ \partial y}-A_j{\partial A_i/ \partial y}=0$$
and
$${\partial \delta_i^{(r)} / \partial x_j}-{\partial \delta_j^{(r)}/
\partial x_i}=0.$$
Therefore we can choose $c^{(r+1)}$ such that
$$\delta_i^{(r)}={\partial c^{(r+1)}/ \partial x_i}$$
for all $i$ and then
$$\{\tilde x_1,\dots ,\tilde x_{i-1},\tilde x_{i+1},\dots ,\tilde x_{n-1},\tilde y\}=
(-1)^{n-i}(A_i^{(0)}+\cdots +A_i^{(r)}+\cdots$$
with $A_i^{(s)}$ polynomial of degree $p$ in $y$ for $s=0,\dots ,r.$

To complete the proof, choose
$\tilde \Omega =d\tilde x_1\wedge \dots\wedge d\tilde x_{n-1}\wedge d\tilde y.$
Then
$\tilde\omega =i_\Lambda\tilde\Omega$ is equal to
$\omega$ multiplied by a function of type
$1+u^{(r+1)}+\cdots$ and we have
$$\tilde\omega =y^p(1+\tilde c^{(r+1)}+\cdots )dy +\sum A_idx_i.$$
We can multiply $\tilde\omega$ by the inverse of
$(1+\tilde c^{(r+1)}+\cdots )$ to get
$$\omega '=y^pdy+\sum A'_idx_i,$$
where $A'_i=A_i^{(0)}+\cdots +A_i^{(r)}+{A'_i}^{(r+1)}+\cdots .$
So, step by step, we obtain the proof of our theorem.
\hfill $\triangle$

In the case $p=2$ the last theorem can be improved.
The integrability condition $\omega_0\wedge d\omega_0=0$
is equivalent to the system of
equations:
$$d\alpha_1=d\alpha_2=0,\ d\alpha_0=\alpha_2\wedge\alpha_1,\
\alpha_0\wedge\alpha_1=\alpha_2\wedge\alpha_0=0.$$
So we can write, at the level of formal series,
$\alpha_1=dg,$ $\alpha_2=dh$ and, since
$d(\alpha_0-h\alpha_1)=0,$
we have $\alpha_0=dk+hdg$ for some function $k$.
Now, the last two equations of our system
give $dk\wedge dg=0$ and $dg\wedge dh=0.$
Using
 Moussu's result (\cite{Moussu}),
 we can conclude that there is a function $f$ whose formal series
 satisfies the relations
$g=a\circ f,$   $h=b\circ f$ and $k=c\circ f,$ for some functions
$a,$ $b$ and $c$ in one variable. So we obtain $$\omega_0=
y^2dy+(c'(f)+b(f)a'(f)+ya'(f)+y^2b'(f))df.$$ This can be
interpreted as follows: $\omega_0$ is the formal pullback of a
2-dimensional 1-form
$y^2dy+(\gamma_0(x)+y\gamma_1(x)+y^2\gamma_2(x))dx$ by a mapping
of the form $(x_1,\dots ,x_{n})\mapsto (f(x_1,\dots
,x_{n-1}),x_n).$

  It seems that $\omega_0$ is a formal
pullback of a 2-dimensional 1-form for any value of $p$.

\section{Quadratic integrable 1-forms}

In this paragraph we will give a classification
of quadratic integrable 1-forms
or, equivalently, a classification of quadratic Nambu tensors of
order $n-1,$ up
to multiplication by a constant.

Let $\Lambda$ be such a quadratic Nambu tensor of order $n-1.$
Its modular tensor $D\Lambda$
relatively to any constant volume form is
intrinsically defined and it is a linear
Nambu tensor of order $n-2.$
The classification of linear Nambu tensors
(\cite{D-N}) says that we have the following two cases.

\noindent 1- {\sl $D\Lambda$ is of type 2:}

This means that we have, in a suitable coordinates system,
$$
D\Lambda =
\partial / \partial x_{4} \wedge \dots \wedge
\partial / \partial x_{n} \wedge X
$$
where $X$ is a vector field depending
on coordinates $x_1,$ $x_2$ and $x_3$ only.

With the notation introduced in definition \ref{2.r}, $D\Lambda$
is of type $2.(n-3)$.
So, due to the generalized Kupka phenomenon
(theorem \ref{KupkaNambu}), $\Lambda$ is
also of type $2.(n-3)$. Then we have
$$\Lambda = {\partial /\partial x_4}
\wedge \dots\wedge{\partial /\partial x_{n}}
\wedge \Lambda_3,$$
where $\Lambda_3$ is a quadratic Poisson structure
depending on the variables
$x_1,$ $x_2$ and $x_3$ only.
We see that the classification of these Nambu structures reduces to the
classification of quadratic 3-dimensional Poisson structures.
The latter classification is known (see \cite{D-H}).

\noindent 2- {\sl $D\Lambda$ is of type 1:}

In this case it is easier to
work with the associated quadratic integrable 1-form $\omega ;$
 $D\Lambda$ is of type 1 if we have $d\omega =dx\wedge dq$
 where $q$ is a quadratic form of type
$q=\sum_{i=1}^{r}\pm y_i^2/2+ xz$ in a system of
coordinates $x,y_1,\dots ,y_r,z,t_1,\dots ,t_s$
with $r+s=n-2$  or $q=\sum_{i=1}^{r}\pm y_i^2/2$
in a system of coordinates
$x,y_1,\dots ,y_r,t_1,\dots ,t_s$ with $r+s=n-1.$
In the sequel we will consider the first
case with $r\geq 2.$ The other cases,
with $r=0,$ $r=1$ or without variable $z$ are easier,
and we let them to the reader.

Since we have  $d\omega =-d(qdx),$ we can express $\omega $ in the form
$\omega = -qdx+df,$ where $f$ is a homogeneous
function of degree 3. Denote $\overline q=\sum_{i=1}^{r}\pm y_i^2/2.$
Then we have
$$0=\omega\wedge d\omega =df\wedge dx\wedge dq=$$$$(
\sum_i{\partial f /\partial {y_i}}dy_i+
{\partial f /\partial z}dz+\sum_j{\partial f /\partial {t_j}}dt_j)\wedge dx\wedge
(d\overline q+ xdz).$$
The terms with $dt_j\wedge dx\wedge dy_i$ in this relation give
$\partial f /\partial {t_j}=0,$ therefore $f$ is
independent of $t_j.$ The terms
with
$dy_j\wedge dx\wedge dy_i$ give
$$\sum_i{\partial f /\partial {y_i}}dy_i\wedge d\overline q=0,$$
and an elementary calculation leads to the relation
$$f=(\lambda x+\mu z)\overline q+ b(x,z).$$
Using this relation we obtain
$$0=((\mu \overline q +{\partial b/\partial z})dz+(\lambda x+\mu z)d\overline q)\wedge dx
\wedge(d\overline q+ xdz)$$$$=((\lambda x+\mu z)x-\mu \overline q-{\partial b/\partial z})
d\overline q\wedge dx\wedge dz.$$
Considering the terms with $y_i$ in the latter relation, we obtain,
step-by-step:
$\mu =0,$
${\partial b/\partial z}=\lambda x^2,$
$b=\lambda x^2z+\alpha x^3$, and finally
$f=\lambda xq+\alpha x^3,$
where $\alpha $ is a constant.

Returning to the expression of $\omega ,$ we get
$$\omega =\theta qdx+\beta xdq+\gamma x^2dx$$
where $\theta$, $\beta$ and $\gamma$ are constants.

The preceding calculations are summarized in the  following theorem.

\begin{theorem} If $\omega$ is a quadratic integrable 1-form,
it is the pull-back
of a 3-dimensional integrable 1-form.
More precisely, if $d\omega$ is of type 2,
then $\omega$ is a  quadratic integrable 1-form depending only on
three (well chosen)
 coordinates; if $d\omega$ is of type 1, then,
 in a suitable system of coordinates,
$$\omega=\phi^*((\gamma x^2+\theta y)dx+\beta xdy)$$
with $$\phi (x_1,\dots ,x_n)=(x_1,\sum_{i=1}^{r}\pm y_i^2/2+\epsilon xz)$$
and $\epsilon$, $\beta$, $\theta$ and $\gamma$
being constants ($\epsilon =0$ or $\epsilon=1$).
In this last case $\omega$ is, in fact,
the pull-back of a 2-dimensional 1-form.
\end{theorem}

\noindent{\bf Conjecture.}
The preceding section and the theorem above lead to the following
conjecture: Every integrable 1-form on ${\mathbb R}^n$
or ${\mathbb C}^n$ with a non-zero
2-jet at 0 is, up to multiplication by a non-vanishing function,
the pull-back of an
integrable 1-form in dimension 3.
More generally we can ask if every integrable 1-form on   ${\mathbb R}^n$ or
${\mathbb C}^n$ with a non-zero
$q$-jet at 0 is, up to multiplication by a non-vanishing function,
the pull-back of an
integrable 1-form in dimension $q+1$.

\bibliographystyle{amsalpha}

D\'epartement de Math\'ematiques, Universit\'e Montpellier II (France)\\
and Department of mathematics, Technion, Haifa (Israel)\\

\end{document}